\documentclass[reqno,11pt]{amsart}

\usepackage[letterpaper,hmargin=1.25in,vmargin=0.9in]{geometry}
\usepackage[T1]{fontenc}
\usepackage{lmodern}
\usepackage{microtype}

\usepackage{mathtools,amssymb}
\DeclareMathOperator{\conv}{conv}
\newtheorem{theorem}{Theorem}[section]
\newtheorem{lemma}[theorem]{Lemma}
\newtheorem{proposition}[theorem]{Proposition}
\newtheorem{conjecture}[theorem]{Conjecture}
\theoremstyle{remark}
\newtheorem{remark}{Remark}[section]
\numberwithin{equation}{section}
\allowdisplaybreaks[1]
\usepackage{xcolor}
\usepackage[unicode]{hyperref}
\hypersetup{
  colorlinks=false,
  pdflinkmargin=1pt,
  citecolor=blue,
  linkcolor=blue,
  urlcolor=blue,
  pdftitle={Difference bodies and lattice covering radii},
  pdfauthor={Yanlu Lian, Jun Wang, and Fei Xue},
  pdfsubject={Mathematics},
  pdfkeywords={lattice covering radius, difference body, integer simplex}
}
\newcommand{\doi}[1]{DOI \href{https://doi.org/#1}{\nolinkurl{#1}}}

\title[Difference bodies and lattice covering radii]{Difference bodies and lattice covering radii}

\author{Yanlu Lian}
\address{School of Mathematics, Hangzhou Normal University, Hangzhou 311121, China}
\email{yllian@hznu.edu.cn}

\author{Jun Wang}
\address{School of Mathematical Sciences \& Institute of Mathematics and Interdisciplinary Sciences, Tianjin Normal University, Tianjin 300387, China}
\email{junwang@tjnu.edu.cn}

\author{Fei Xue}
\address{School of Mathematical Sciences, Nanjing Normal University, Nanjing 210046, China}
\email{05429@njnu.edu.cn}

\date{}

\begin{document}
\begin{abstract}
We construct integer simplices $P_h$ in dimensions $n=2^h$ and compute the covering radii of $P_h$ and $P_h-P_h$ with respect to $\mathbb Z^n$ exactly.
Along these dimensions, their ratio is asymptotic to $\sqrt{(2\pi/\log 2)\log n}$.
Taking products with cubes then shows that Dadush's symmetrization parameter tends to infinity with the dimension, answering the question recorded by Reis and Rothvoss.
We also prove a logarithmic upper bound for this ratio for arbitrary convex bodies and full-rank lattices, using an exponential covering estimate and the Ruzsa triangle inequality on a compact abelian group.
\end{abstract}
\maketitle

\section{Introduction}

Let $K\subset\mathbb R^n$ be a convex body and let $\Lambda\subset\mathbb R^n$ be a full-rank lattice.
The covering radius of $K$ with respect to $\Lambda$ is
$$\mu(K,\Lambda)=\min\{t>0:\Lambda+tK=\mathbb R^n\}.$$
The difference body $K-K=\{x-y:x,y\in K\}$ is symmetric about the origin.
Replacing $K$ by $K-K$ cannot increase the covering radius and decreases it by a factor of two when $K$ is centrally symmetric. We study how large this factor can be for general convex bodies. For background on covering radii and their place in the geometry of numbers, we refer to \cite{GL87,AH23,Rot26}. Recent work includes results on covering radii of three-dimensional lattice zonotopes \cite{ACS26} and Brunn--Minkowski type inequalities for covering radii of Minkowski sums in the plane \cite{CFK26}.

Following Dadush \cite{Dad19}, we consider the symmetrization parameter
\begin{equation}
\label{eq:symmetrization-constant}
 C_{\mathrm{sym}}(n)
 =\sup_{K,\Lambda}
   \frac{\mu(K,\Lambda)}{\mu(K-K,\Lambda)},
\end{equation}
where the supremum is over all convex bodies and full-rank lattices in $\mathbb R^n$. By an invertible linear transformation, we may restrict the supremum to $\Lambda=\mathbb Z^n$. We use natural logarithms, with $\log_2$ reserved for base two.

Dadush \cite{Dad19} raised the possibility that $C_{\mathrm{sym}}(n)$ might be bounded by an absolute constant.
Reis and Rothvoss \cite{RR23} observed that the standard simplex gives $C_{\mathrm{sym}}(n)\geq 4n/(n+1)$.
They also noted that it was unknown whether $C_{\mathrm{sym}}(n)$ tends to infinity with $n$.

Upper bounds for $C_{\mathrm{sym}}(n)$ have been obtained through the projection lower bounds of Kannan and Lov\'asz \cite{KL88}, which are determined by the volumes of projected bodies and the determinants of projected lattices. The projections are taken onto nonzero subspaces for which the projected lattice is discrete. By the Rogers--Shephard inequality \cite{RS57}, replacing $K$ by $K-K$ decreases each such lower bound by at most a factor of four. Consequently, a uniform factor comparing the covering radius with the largest projection lower bound also bounds $C_{\mathrm{sym}}(n)$, up to this factor of four \cite{Dad19}. The comparison of Reis and Rothvoss \cite{RR23} implies $C_{\mathrm{sym}}(n)=O(\log^3(2n))$, as also recorded in \cite[Theorem~5.6]{Rot26}. The sharper comparison of Letwin, Reis and Rothvoss \cite[Corollary~1.3]{LRR26} yields $C_{\mathrm{sym}}(n)=O(\log^{3/2}(2n))$.

We show that $C_{\mathrm{sym}}(n)$ is unbounded by constructing integer simplices with the following asymptotic ratio.

\begin{theorem}
\label{thm:main}
There is a sequence of integer simplices $P_h\subset\mathbb R^n$, with $n=2^h$ and $h\geq0$, such that, as $h\to\infty$,
\begin{equation}
\label{eq:asymptotic-ratio}
 \frac{\mu(P_h,\mathbb Z^n)}{\mu(P_h-P_h,\mathbb Z^n)}
 \sim\sqrt{2\pi h}
 =\sqrt{\frac{2\pi}{\log2}\log n}.
\end{equation}
\end{theorem}

The exact values of both covering radii are given in \eqref{eq:exact-radii}.
Cartesian products with suitable cubes preserve the two radii and extend the lower bound to every dimension, giving $C_{\mathrm{sym}}(n)\to\infty$ and answering the question of Reis and Rothvoss.

We also prove a logarithmic upper bound, improving the $O(\log^{3/2}(2n))$ estimate obtained from the projection comparison.

\begin{proposition}
\label{prop:log-upper}
For every $n\geq1$, every convex body $K\subset\mathbb R^n$ and full-rank lattice $\Lambda\subset\mathbb R^n$,
\begin{equation}
\label{eq:log-upper}
 \mu(K,\Lambda)
 \leq\bigl(2\lfloor\log_2 n\rfloor+10\bigr)
       \mu(K-K,\Lambda).
\end{equation}
In particular, $C_{\mathrm{sym}}(n)=O(\log(2n))$.
\end{proposition}

The parameter $C_{\mathrm{sym}}(n)$ also appears in the running-time bound for Dadush's randomized integer programming algorithm \cite[Theorem~2.10]{Dad19}, which, with high probability, finds an integer point in an appropriately represented convex body or correctly determines that none exists. Proposition~\ref{prop:log-upper} improves the bound on the symmetrization parameter in that estimate to $O(\log(2n))$.

The paper is organised as follows. Section~\ref{sec:construction} gives the three-dimensional example, the recursive construction and the proof of Theorem~\ref{thm:main}. Section~\ref{sec:log-upper} proves Proposition~\ref{prop:log-upper}.
Section~\ref{sec:discussion} extends the examples to all dimensions and formulates a conjecture on the asymptotic order of $C_{\mathrm{sym}}(n)$.

\section{The recursive simplices and their covering radii}
\label{sec:construction}

By the Minkowski--Radon bound on asymmetry, every convex body $K\subset\mathbb R^n$ can be translated so that $-K\subseteq nK$ (see \cite{Rad16} and \cite[Lemma~2.2]{ABG18}).
For this translate we have $K-K\subseteq(n+1)K$.
Since translation preserves the covering radius and enlarging a body can only decrease it, we obtain
$$
\mu(K,\Lambda)\leq(n+1)\mu(K-K,\Lambda).
$$
The tetrahedron in Section~\ref{subsec:tetrahedron} attains this upper bound in dimension three. For sufficiently large $n$, however, the logarithmic upper bound established in Section~\ref{sec:log-upper} gives $C_{\mathrm{sym}}(n)<n+1$, so the bound $n+1$ is not sharp in those dimensions.

The following lemma reduces both covering-radius computations to finding an axis-parallel unit cube whose vertices lie on the boundary of a suitable dilate.

\begin{lemma}
\label{lem:unit-cube}
Let $K\subset\mathbb R^n$ be a convex body, let $r>0$, and let $\boldsymbol a\in\mathbb R^n$. If every vertex of the unit cube
$$
C=\boldsymbol a+[-1/2,1/2]^n
$$
lies on the boundary of $rK$, then
$$
\mu(K,\mathbb Z^n)=r.
$$
\end{lemma}

\begin{proof}
Since $C$ is the convex hull of its vertices and $rK$ is convex, we have $C\subseteq rK$. The integer translates of $C$ cover the whole space, so
$$
rK+\mathbb Z^n=\mathbb R^n.
$$

Fix $0<s<r$ and suppose that $sK+\mathbb Z^n=\mathbb R^n$.
Since $C\subseteq rK$, the centre $\boldsymbol a$ of $C$ lies in the interior of $rK$.
The open segment joining $\boldsymbol a$ to any other point of $rK$ lies in its interior, so $0<s/r<1$ gives
$$
\boldsymbol a+\frac sr(rK-\boldsymbol a)
=sK+\left(1-\frac sr\right)\boldsymbol a
\subseteq\operatorname{int}(rK).
$$
Translating the assumed covering by $\left(1-s/r\right)\boldsymbol a$ still gives a covering of the whole space. In particular, there is a vector $\boldsymbol k\in\mathbb Z^n$ such that
$$
\boldsymbol a+\frac12\mathbf1
\in sK+\left(1-\frac sr\right)\boldsymbol a+\boldsymbol k,
$$
where $\mathbf1=(1,\ldots,1)\in\mathbb R^n$.
Thus
$$
\boldsymbol y:=\boldsymbol a+\frac12\mathbf1-\boldsymbol k
\in\operatorname{int}(rK).
$$
For each coordinate, $y_i-a_i=1/2-k_i\in1/2+\mathbb Z$, and hence $|y_i-a_i|\geq1/2$. Choose the vertex $\boldsymbol v$ of $C$ on the same side of $\boldsymbol a$ as $\boldsymbol y$ in each coordinate, so that
$$
v_i=
\begin{cases}
a_i+1/2,&y_i>a_i,\\
a_i-1/2,&y_i<a_i.
\end{cases}
$$
If $\boldsymbol y=\boldsymbol v$, then $\boldsymbol v\in\operatorname{int}(rK)$, a contradiction. We may therefore assume $\boldsymbol y\neq\boldsymbol v$.

Extend the segment from $\boldsymbol y$ to $\boldsymbol v$ slightly beyond $\boldsymbol v$. Each coordinate with $|y_i-a_i|>1/2$ moves into the interior of $[a_i-1/2,a_i+1/2]$ after passing $v_i$, while each coordinate with $|y_i-a_i|=1/2$ stays fixed. We can therefore choose a point $\boldsymbol z\in C$ on the extension such that $\boldsymbol v$ lies strictly between $\boldsymbol y$ and $\boldsymbol z$. Since $\boldsymbol y\in\operatorname{int}(rK)$ and $\boldsymbol z\in C\subseteq rK$, convexity gives $\boldsymbol v\in\operatorname{int}(rK)$.
This contradicts the assumption that $\boldsymbol v$ lies on the boundary of $rK$. Thus no positive scale smaller than $r$ gives a covering, and we conclude that $\mu(K,\mathbb Z^n)=r$.
\end{proof}

\subsection{A three-dimensional equality example}
\label{subsec:tetrahedron}

Consider the tetrahedron
$$
Q=\conv\{(1,2,0),(1,-2,0),(-1,0,2),(-1,0,-2)\}.
$$

\begin{proposition}
\label{prop:tetrahedron}
The tetrahedron $Q$ satisfies
$$
\mu(Q,\mathbb Z^3)=1,
\qquad
\mu(Q-Q,\mathbb Z^3)=\frac14.
$$
Consequently, $C_{\mathrm{sym}}(3)=4$.
\end{proposition}

\begin{proof}
Put $C=[-1/2,1/2]^3$. We first use the sections of $Q$ to show that every vertex of $C$ lies on the boundary of $Q$.

Write $Q$ as the convex hull of two segments,
$$
Q=\conv(\ell_1\cup\ell_2),
$$
where
$$
\ell_1=\{1\}\times[-2,2]\times\{0\},
\qquad
\ell_2=\{-1\}\times\{0\}\times[-2,2].
$$
Since each segment is convex, every point $\boldsymbol x\in Q$ is a convex combination of one point from each segment. Thus there are points
$$
\boldsymbol p=(1,u,0)\in\ell_1,
\qquad
\boldsymbol q=(-1,0,v)\in\ell_2,
\qquad
-2\leq u,v\leq2,
$$
and a number $\lambda\in[0,1]$ such that
$$
\boldsymbol x
=\lambda\boldsymbol p+(1-\lambda)\boldsymbol q
=(2\lambda-1,\lambda u,(1-\lambda)v).
$$
If the first coordinate of $\boldsymbol x$ is $a$, then $2\lambda-1=a$, so
$$
\lambda=\frac{1+a}{2}.
$$
Since $u$ and $v$ vary independently over $[-2,2]$, we obtain
\begin{equation}
\label{eq:tetrahedron-section}
Q\cap\{x_1=a\}
=\{a\}\times[-(1+a),1+a]\times[-(1-a),1-a]
\end{equation}
for $-1\leq a\leq1$.
Taking $a=1/2$ and $a=-1/2$ in \eqref{eq:tetrahedron-section} shows that all eight vertices of $C$ belong to $Q$.
The four vertices with first coordinate $1/2$ satisfy $|x_3|=1-x_1$, and the other four satisfy $|x_2|=1+x_1$. Since every point of $Q$ satisfies
$$
|x_3|\leq1-x_1,
\qquad
|x_2|\leq1+x_1,
$$
each of these vertices lies on the boundary of $Q$. Lemma~\ref{lem:unit-cube} therefore gives
$$
\mu(Q,\mathbb Z^3)=1.
$$

Now put $D=Q-Q$. The first coordinate of every point of $Q$ belongs to $[-1,1]$, so the first coordinate of the difference of two points of $Q$ is at most $2$. It equals $2$ if and only if the first point belongs to $\ell_1$ and the second belongs to $\ell_2$. Consequently,
$$
D\cap\{x_1=2\}
=\ell_1-\ell_2
=\{2\}\times[-2,2]^2.
$$
Since $D$ is contained in the half-space $x_1\leq2$, this square lies on its boundary. By the symmetry $D=-D$, the opposite square $\{-2\}\times[-2,2]^2$ also lies on the boundary of $D$. Scaling the eight vertices of these squares by $1/4$ gives the eight vertices of $C$. Thus every vertex of $C$ lies on the boundary of $\frac14D$, and Lemma~\ref{lem:unit-cube} gives
$$
\mu(Q-Q,\mathbb Z^3)=\frac14.
$$

The ratio for $Q$ is $4$, and the elementary bound $C_{\mathrm{sym}}(3)\leq4$ therefore gives $C_{\mathrm{sym}}(3)=4$.
\end{proof}

\begin{remark}
\label{rem:tetrahedron}
Although Proposition~\ref{prop:tetrahedron} also follows from the arguments below, we retain this proof to explain the construction and the argument used in higher dimensions.
\end{remark}

As in the preceding proof, the tetrahedron $Q$ can also be written as
$$
Q=\operatorname{conv}\Bigl((\{1\}\times2[-1,1]\times\{0\})\cup(\{-1\}\times\{0\}\times2[-1,1])\Bigr).
$$
Replacing the interval $[-1,1]$ by a simplex in this representation leads to the recursive construction below.

\subsection{The recursive construction}
\label{subsec:recursive-construction}

Put $m_h=2^h-1$, let $C_0$ be the unique point of $\mathbb R^0$, and define recursively
\begin{equation}
\label{eq:recursive-base}
C_h=\conv\Bigl(
(\{1\}\times2C_{h-1}\times\{\boldsymbol0\})
\cup
(\{-1\}\times\{\boldsymbol0\}\times2C_{h-1})
\Bigr),
\end{equation}
where the zero vectors belong to $\mathbb R^{m_{h-1}}$. Thus $C_h\subset\mathbb R^{m_h}$, with $C_1=[-1,1]$ and $C_2=Q$.
Every $C_h$ contains the origin. For if $\boldsymbol0\in C_{h-1}$, then the two points
$$
(1,\boldsymbol0,\boldsymbol0),
\qquad
(-1,\boldsymbol0,\boldsymbol0)
$$
belong to $C_h$, and their midpoint is the origin.

Set
\begin{equation}
\label{eq:recursive-simplex}
P_h=\conv\bigl(\{\boldsymbol0\}\cup(\{1\}\times C_h)\bigr)
\subset\mathbb R^{2^h},
\qquad
D_h=P_h-P_h.
\end{equation}
We write $n=2^h=m_h+1$.

For $r>0$, every point of $rC_h$ can be written as
$$
\lambda(r,2r\boldsymbol u,\boldsymbol0)
+(1-\lambda)(-r,\boldsymbol0,2r\boldsymbol v),
$$
where $0\leq\lambda\leq1$ and $\boldsymbol u,\boldsymbol v\in C_{h-1}$.
Its first coordinate is $a=r(2\lambda-1)$, so
$$
\lambda=\frac{r+a}{2r}.
$$
Substitution gives
$$
\bigl(a,(r+a)\boldsymbol u,(r-a)\boldsymbol v\bigr).
$$
Since $\boldsymbol u$ and $\boldsymbol v$ vary independently over $C_{h-1}$,
\begin{equation}\label{eq:base-section}
rC_h\cap\{x_1=a\}
=\{a\}\times(r+a)C_{h-1}\times(r-a)C_{h-1},
\qquad |a|\leq r.
\end{equation}
The section is empty when $|a|>r$.

We verify by induction that $C_h$ is an integer simplex of dimension $m_h$.
This holds for the point $C_0$.
Suppose that $C_{h-1}$ is an integer simplex of dimension $m_{h-1}$.
Taking $r=1$ and $a=0$ in \eqref{eq:base-section} gives the central section $\{0\}\times C_{h-1}\times C_{h-1}$, which has dimension $2m_{h-1}$.
The point $(1,\boldsymbol0,\boldsymbol0)$ lies outside the affine hull of this section, so $\dim C_h=1+2m_{h-1}=m_h$.
Since $C_h$ is the convex hull of $m_h+1$ integer points, it is an integer simplex.
Finally, the origin lies outside the affine hull of $\{1\}\times C_h$, so $P_h$ is an $n$-dimensional integer simplex.

\subsection{Vertex coefficients and difference bodies}
\label{subsec:vertex-coefficients}

We first describe the difference body of a simplex in terms of linear combinations of its vertices. We then compute the vertex coefficients needed to apply Lemma~\ref{lem:unit-cube} to $P_h$ and $D_h$.

\begin{lemma}
\label{lem:simplex-difference}
Let $P=\conv\{\boldsymbol p_0,\ldots,\boldsymbol p_n\}$ be an $n$-dimensional simplex in $\mathbb R^n$. Then the following assertions hold.
\begin{enumerate}
\item
For every $\boldsymbol x\in\mathbb R^n$, there are unique real coefficients $\lambda_0,\ldots,\lambda_n$ satisfying
$$
\boldsymbol x=\sum_{i=0}^n\lambda_i\boldsymbol p_i,
\qquad
\sum_{i=0}^n\lambda_i=0.
$$
\item
For these coefficients and every $r\geq0$,
\begin{equation}
\label{eq:simplex-difference-coefficients}
\boldsymbol x\in r(P-P)
\quad\Longleftrightarrow\quad
\frac12\sum_{i=0}^n|\lambda_i|\leq r.
\end{equation}
\end{enumerate}
\end{lemma}

\begin{proof}
Since $P$ is an $n$-dimensional simplex, the vectors $\boldsymbol p_i-\boldsymbol p_0$, $1\leq i\leq n$, form a basis of $\mathbb R^n$.
Thus $\boldsymbol x$ has a unique representation
$$
\boldsymbol x=\sum_{i=1}^n\lambda_i(\boldsymbol p_i-\boldsymbol p_0).
$$
Putting $\lambda_0=-\sum_{i=1}^n\lambda_i$ proves the first assertion.

If $\boldsymbol x\in r(P-P)$, the definition of the difference body gives points $\boldsymbol u,\boldsymbol v\in P$ such that
$$
\boldsymbol x=r(\boldsymbol u-\boldsymbol v).
$$
Since $\boldsymbol u,\boldsymbol v\in P$, there are nonnegative coefficients $\alpha_i,\beta_i$ satisfying
$$
\boldsymbol u=\sum_{i=0}^n\alpha_i\boldsymbol p_i,
\qquad
\boldsymbol v=\sum_{i=0}^n\beta_i\boldsymbol p_i,
\qquad
\sum_{i=0}^n\alpha_i=\sum_{i=0}^n\beta_i=1.
$$
Substitution gives
$$
\boldsymbol x=\sum_{i=0}^n r(\alpha_i-\beta_i)\boldsymbol p_i.
$$
These coefficients sum to zero, so uniqueness implies $\lambda_i=r(\alpha_i-\beta_i)$. Using $|\alpha_i-\beta_i|\leq\alpha_i+\beta_i$, we obtain
$$
\frac12\sum_{i=0}^n|\lambda_i|
=\frac r2\sum_{i=0}^n|\alpha_i-\beta_i|
\leq\frac r2\sum_{i=0}^n(\alpha_i+\beta_i)
=r.
$$

Conversely, suppose that $\frac12\sum_i|\lambda_i|\leq r$. Since $\sum_i\lambda_i=0$, the sum of the positive coefficients equals the sum of the absolute values of the negative coefficients. Put
$$
q=\sum_{\lambda_i>0}\lambda_i
=\sum_{\lambda_i<0}(-\lambda_i)
=\frac12\sum_{i=0}^n|\lambda_i|
\leq r.
$$
If $q=0$, then $\boldsymbol x=\boldsymbol0$ and the assertion follows. If $q>0$, then
$$
\boldsymbol x
=q\left(
\sum_{\lambda_i>0}\frac{\lambda_i}{q}\boldsymbol p_i
-
\sum_{\lambda_i<0}\frac{-\lambda_i}{q}\boldsymbol p_i
\right).
$$
Both sums in parentheses are convex combinations of vertices of $P$, so $\boldsymbol x\in q(P-P)$. Since $P-P$ is convex and contains the origin, $q\leq r$ implies $q(P-P)\subseteq r(P-P)$. This proves the reverse implication in \eqref{eq:simplex-difference-coefficients}.
\end{proof}

\begin{lemma}
\label{lem:vertex-coefficients}
For $h\geq1$, let $\boldsymbol p_1,\ldots,\boldsymbol p_n$ be the nonzero vertices of $P_h$, where $n=2^h$. For every $t\in\mathbb R$ and $\boldsymbol w\in\{-1/2,1/2\}^{m_h}$, consider the unique representation
$$
(t,\boldsymbol w)=\sum_{i=1}^n\gamma_i\boldsymbol p_i.
$$
Up to reordering, the coefficients $\gamma_1,\ldots,\gamma_n$ are
\begin{equation}
\label{eq:vertex-coefficients}
2^{-h}\left(t+\frac12\sum_{j=1}^h\varepsilon_j\right),
\qquad
(\varepsilon_1,\ldots,\varepsilon_h)\in\{-1,1\}^h,
\end{equation}
with each sign sequence occurring once.
\end{lemma}

\begin{proof}
Existence and uniqueness follow from Lemma~\ref{lem:simplex-difference}, since the remaining vertex of $P_h$ is the origin. We prove the coefficient formula by induction on $h$. For $h=1$, the two nonzero vertices are $(1,1)$ and $(1,-1)$. For $(t,a)$ with $a=\pm1/2$, comparison of the coordinates gives the coefficients
$$
\frac{t+a}{2},
\qquad
\frac{t-a}{2}.
$$
These are the two values in \eqref{eq:vertex-coefficients} for $h=1$.

Let $h\geq2$ and assume the assertion for $h-1$. If $\boldsymbol c_1,\ldots,\boldsymbol c_{2^{h-1}}$ are the vertices of $C_{h-1}$, the nonzero vertices of $P_h$ are
$$
(1,1,2\boldsymbol c_i,\boldsymbol0),
\qquad
(1,-1,\boldsymbol0,2\boldsymbol c_i),
\qquad
1\leq i\leq2^{h-1}.
$$
Since $m_h=1+2m_{h-1}$, we divide $\boldsymbol w=(w_1,\ldots,w_{m_h})$ into $(a,\boldsymbol u,\boldsymbol v)$ according to the positions of $2\boldsymbol c_i$ in these two groups, where
$$
\begin{aligned}
a&=w_1,\\
\boldsymbol u&=(w_2,\ldots,w_{m_{h-1}+1}),\\
\boldsymbol v&=(w_{m_{h-1}+2},\ldots,w_{m_h}).
\end{aligned}
$$
Each of $\boldsymbol u$ and $\boldsymbol v$ has $m_{h-1}$ coordinates, and
$$
a\in\{-1/2,1/2\},
\qquad
\boldsymbol u,\boldsymbol v\in\{-1/2,1/2\}^{m_{h-1}}.
$$

Let $\alpha_i$ and $\beta_i$ be the coefficients corresponding to the two groups of vertices, so that
$$
(t,a,\boldsymbol u,\boldsymbol v)
=\sum_i\alpha_i(1,1,2\boldsymbol c_i,\boldsymbol0)
+\sum_i\beta_i(1,-1,\boldsymbol0,2\boldsymbol c_i).
$$
Comparing coordinates, this equality is equivalent to the system
$$
\begin{aligned}
\sum_i(\alpha_i+\beta_i)&=t,
&\qquad \sum_i(\alpha_i-\beta_i)&=a,\\
\sum_i2\alpha_i\boldsymbol c_i&=\boldsymbol u,
&\qquad \sum_i2\beta_i\boldsymbol c_i&=\boldsymbol v.
\end{aligned}
$$
Adding and subtracting the first two equations shows that this system is equivalent to
$$
(t+a,\boldsymbol u)=\sum_i2\alpha_i(1,\boldsymbol c_i),
\qquad
(t-a,\boldsymbol v)=\sum_i2\beta_i(1,\boldsymbol c_i).
$$
The vectors $(1,\boldsymbol c_i)$ are the nonzero vertices of $P_{h-1}$, so the inductive hypothesis applies to both representations.

Dividing the coefficients given by the inductive hypothesis by $2$ gives the two lists
$$
2^{-h}\left(t+a+\frac12\sum_{j=1}^{h-1}\varepsilon_j\right),
\qquad
2^{-h}\left(t-a+\frac12\sum_{j=1}^{h-1}\varepsilon_j\right),
$$
each ranging over all $(\varepsilon_1,\ldots,\varepsilon_{h-1})\in\{-1,1\}^{h-1}$. Since $a=\pm1/2$, the two lists together correspond to $\varepsilon_h=\pm1$ in \eqref{eq:vertex-coefficients}, with each sign sequence occurring once.
This completes the induction.
\end{proof}

\subsection{Proof of Theorem~\ref{thm:main}}

\begin{proof}
For $h=0$, we have $P_0=[0,1]$ and $D_0=[-1,1]$, whose covering radii are $1$ and $1/2$, respectively. We now assume $h\geq1$.

To compute the covering radius of $P_h$, put $r_h=h/2+1$. Since $P_h$ is an $n$-dimensional simplex with the origin as one vertex, every point $\boldsymbol x=(t,\boldsymbol w)$ has a unique representation
$$
\boldsymbol x=\sum_{i=1}^n\gamma_i\boldsymbol p_i
$$
in terms of its nonzero vertices. Each nonzero vertex has first coordinate $1$, so $\sum_i\gamma_i=t$. To decide whether $\boldsymbol x$ belongs to $r_hP_h$, it is enough to write $\boldsymbol x/r_h$ as a convex combination of the vertices of $P_h$. We have
$$
\frac{\boldsymbol x}{r_h}
=
\sum_{i=1}^n\frac{\gamma_i}{r_h}\boldsymbol p_i
+
\left(1-\frac{t}{r_h}\right)\boldsymbol0.
$$
These coefficients sum to $1$, so $\boldsymbol x\in r_hP_h$ if and only if $\gamma_i\geq0$ for all $i$ and $t\leq r_h$.

Consider the unit cube
$$
U_h=\left[\frac h2,\frac h2+1\right]\times[-1/2,1/2]^{m_h}.
$$
Its centre is $\bigl((h+1)/2,\boldsymbol0\bigr)$. Its vertices satisfy $t\in\{h/2,h/2+1\}$ and $\boldsymbol w\in\{-1/2,1/2\}^{m_h}$. By Lemma~\ref{lem:vertex-coefficients}, the smallest coefficient at each such vertex is
$$
\min_{1\leq i\leq n}\gamma_i
=2^{-h}\left(t-\frac h2\right).
$$
All coefficients are therefore nonnegative and their sum is $t\leq r_h$, so every vertex of $U_h$ belongs to $r_hP_h$.
When $t=h/2$, the smallest coefficient is zero, and the point lies on the face opposite the corresponding nonzero vertex of $r_hP_h$. When $t=h/2+1=r_h$, the weight $1-t/r_h$ at the origin is zero, and the point lies on the face opposite the origin. Since $r_hP_h$ is an $n$-dimensional simplex, these faces lie on its boundary. Thus every vertex of $U_h$ lies on the boundary of $r_hP_h$, and Lemma~\ref{lem:unit-cube} gives
\begin{equation}
\label{eq:simplex-radius}
\mu(P_h,\mathbb Z^n)=r_h=\frac{h+2}{2}.
\end{equation}

We now determine the covering radius of $D_h$.
Let $S_h=\varepsilon_1+\cdots+\varepsilon_h$, where the $\varepsilon_j$ are independent signs taking $1$ and $-1$ with equal probability, and put $S_0=0$. The symbol $\mathbb E$ denotes expectation.

Take any vertex
$$
\boldsymbol v
=\left(\frac{\delta_0}{2},\ldots,\frac{\delta_{n-1}}2\right),
\qquad
\delta_i\in\{-1,1\},
$$
of the unit cube $[-1/2,1/2]^n$.
Apply Lemma~\ref{lem:simplex-difference} to the vertices $\boldsymbol p_0=\boldsymbol0,\boldsymbol p_1,\ldots,\boldsymbol p_n$ of $P_h$.
The coefficients associated with $\boldsymbol v$ satisfy
$$
\boldsymbol v=\sum_{i=0}^n\lambda_i\boldsymbol p_i,
\qquad
\sum_{i=0}^n\lambda_i=0.
$$
Every nonzero vertex has first coordinate $1$. Comparing first coordinates and using the zero sum of the coefficients gives
$$
\sum_{i=1}^n\lambda_i=\frac{\delta_0}{2},
\qquad
\lambda_0=-\frac{\delta_0}{2}.
$$
Since $\boldsymbol p_0=\boldsymbol0$, the remaining coefficients $\lambda_1,\ldots,\lambda_n$ are those in the representation with respect to the nonzero vertices from Lemma~\ref{lem:vertex-coefficients}.
Taking $t=\delta_0/2$ in \eqref{eq:vertex-coefficients}, we see that, up to reordering, they are
\begin{equation}
\label{eq:cube-vertex-coefficients}
2^{-h-1}\left(\delta_0+\sum_{j=1}^h\varepsilon_j\right),
\qquad
(\varepsilon_1,\ldots,\varepsilon_h)\in\{-1,1\}^h.
\end{equation}
When $\delta_0=-1$, replacing every $\varepsilon_j$ by its negative shows that the sum of the absolute values in \eqref{eq:cube-vertex-coefficients} is the same as for $\delta_0=1$. Since $|\lambda_0|=1/2$ and each sign sequence has probability $2^{-h}$, we obtain
\begin{equation}
\label{eq:cube-vertex-scale}
\begin{aligned}
\frac12\sum_{i=0}^n|\lambda_i|
&=\frac14+2^{-h-2}\sum_{\boldsymbol\varepsilon\in\{-1,1\}^h}
\left|1+\sum_{j=1}^h\varepsilon_j\right|\\
&=\frac{1+\mathbb E|1+S_h|}{4}.
\end{aligned}
\end{equation}
Denote this value, which is independent of $\boldsymbol v$, by $\rho_h$. By \eqref{eq:simplex-difference-coefficients}, for every cube vertex $\boldsymbol v$ and every $r\geq0$,
\begin{equation}
\label{eq:cube-vertex-threshold}
\boldsymbol v\in rD_h
\quad\Longleftrightarrow\quad
r\geq\rho_h.
\end{equation}

It follows from \eqref{eq:cube-vertex-threshold} that all vertices of the cube belong to $\rho_hD_h$. If a vertex $\boldsymbol v$ were in the interior of $\rho_hD_h$, then $(1+\eta)\boldsymbol v$ would still belong to $\rho_hD_h$ for sufficiently small $\eta>0$. This would give
$$
\boldsymbol v\in\frac{\rho_h}{1+\eta}D_h,
$$
contrary to \eqref{eq:cube-vertex-threshold}. Thus every vertex lies on the boundary of $\rho_hD_h$, and Lemma~\ref{lem:unit-cube} gives $\mu(D_h,\mathbb Z^n)=\rho_h$. We have therefore proved, for every $h\geq0$, that
\begin{equation}
\label{eq:exact-radii}
\mu(P_h,\mathbb Z^n)=\frac{h+2}{2},
\qquad
\mu(P_h-P_h,\mathbb Z^n)=\rho_h=\frac{1+\mathbb E|1+S_h|}{4}.
\end{equation}

It remains to estimate the ratio. Symmetry of $S_h$ and averaging over the independent sign $\varepsilon_{h+1}$ give
$$
\mathbb E|S_{h+1}|
=\frac{\mathbb E|S_h+1|+\mathbb E|S_h-1|}{2}
=\mathbb E|1+S_h|.
$$
If exactly $k$ of the $h+1$ signs are positive, then $S_{h+1}=2k-h-1$, an event of probability $2^{-(h+1)}\binom{h+1}{k}$.
By symmetry of $S_{h+1}$, we can sum over its positive values and multiply by $2$. The identity
$$
(2k-h-1)\binom{h+1}{k}
=(h+1)\left(\binom h{k-1}-\binom hk\right)
$$
makes this sum telescope and gives
$$
\begin{aligned}
\mathbb E|1+S_h|
&=2^{-h}\sum_{k=\lfloor(h+1)/2\rfloor+1}^{h+1}
(2k-h-1)\binom{h+1}{k}\\
&=(h+1)2^{-h}\binom h{\lfloor h/2\rfloor},
\end{aligned}
$$
where we also used $\binom hk=\binom h{h-k}$. This also follows from the formula for the binomial mean absolute deviation recalled in \cite{BK13}. Stirling's formula now yields
$$
\mathbb E|1+S_h|\sim\sqrt{\frac{2h}{\pi}}.
$$
Combining this with \eqref{eq:exact-radii} and $n=2^h$, we conclude that
$$
\frac{\mu(P_h,\mathbb Z^n)}{\mu(P_h-P_h,\mathbb Z^n)}
=\frac{2(h+2)}{1+\mathbb E|1+S_h|}
\sim\sqrt{2\pi h}
=\sqrt{\frac{2\pi}{\log2}\log n}.
$$
This proves Theorem~\ref{thm:main}.
\end{proof}

\section{A logarithmic upper bound}
\label{sec:log-upper}

The lattice translates of a set cover $\mathbb R^n$ if and only if its image in $\mathbb R^n/\Lambda$ is the whole quotient.
We first record a geometric covering estimate and an inequality for sums of sets in a compact group.

For convex bodies $E,F\subset\mathbb R^n$, let $N(E,F)$ be the least number of translates of $F$ needed to cover $E$.
We use the covering estimate from the proof of Theorem~4.3.7 in \cite{Dad12},
\begin{equation}
\label{eq:exponential-cover}
 N(K-K,K)\leq24^n.
\end{equation}

The first inequality in the next lemma is a special case of the Ruzsa triangle inequality for compact abelian groups \cite[Lemma~2.8]{CGdR19}. We give a direct proof and then derive the estimate for repeated sums that will be used below. A related iteration in the finite setting appears in the proof of \cite[Proposition~3.4]{Skr24}.

\begin{lemma}
\label{lem:torus-growth}
Let $G$ be a compact abelian group with Haar probability measure $m$, and let $A,B\subset G$ be compact sets with $A-A=G$. Then
\begin{equation}
\label{eq:group-growth}
 m(A+B)^2\geq m(B).
\end{equation}
If $a=m(A)>0$, then, for every integer $j\geq1$,
\begin{equation}
\label{eq:iterated-growth}
 m(jA)\geq a^{1/2^{j-1}},
\end{equation}
where $jA$ denotes the sum of $j$ copies of $A$.
\end{lemma}

\begin{proof}
Put $U=A+B$.
For each $d\in G$, choose $a_1,a_2\in A$ with $d=a_1-a_2$. For every $b\in B$, we have
$$a_2+b\in U, \qquad (a_2+b)+d=a_1+b\in U.$$
Hence
$$a_2+B\subset U\cap(U-d), \qquad m\bigl(U\cap(U-d)\bigr)\geq m(B).$$
Write $\mathbf1_U$ for the indicator function of $U$. Fubini's theorem and translation invariance of $m$ give
$$
 \begin{aligned}
 \int_G m\bigl(U\cap(U-d)\bigr)\,dm(d)
 &=\int_G\int_G\mathbf1_U(x)\mathbf1_U(x+d)\,dm(x)\,dm(d)\\
 &=\int_G\mathbf1_U(x)
       \left(\int_G\mathbf1_U(x+d)\,dm(d)\right)dm(x)\\
 &=m(U)\int_G\mathbf1_U(x)\,dm(x)=m(U)^2.
 \end{aligned}
$$
Integrating the preceding lower bound and using $m(G)=1$ therefore gives $m(U)^2\geq m(B)$, which is \eqref{eq:group-growth}.

For $j\geq2$, take $B=(j-1)A$ in \eqref{eq:group-growth} to obtain
$$m(jA)\geq\sqrt{m((j-1)A)}.$$
The estimate \eqref{eq:iterated-growth} follows by induction from $m(A)=a$.
\end{proof}

\begin{proof}[Proof of Proposition~\ref{prop:log-upper}]
Put $D=K-K$ and $\lambda=\mu(D,\Lambda)$. By the definition of the covering radius, $\lambda D+\Lambda=\mathbb R^n$.

Let $\pi:\mathbb R^n\to G=\mathbb R^n/\Lambda$ be the quotient map. Equip $G$ with its Haar probability measure $m$ and put $A=\pi(\lambda K)$. Thus $m(A)$ is the proportion of a fundamental parallelepiped of $\Lambda$ covered by $\lambda K+\Lambda$. The covering identity and the additivity of $\pi$ give
$$A-A=\pi(\lambda D)=G.$$

By \eqref{eq:exponential-cover}, there are vectors $\boldsymbol t_1,\ldots,\boldsymbol t_M\in\mathbb R^n$ with $M\leq24^n$ such that
$$ D\subseteq\bigcup_{i=1}^M(\boldsymbol t_i+K).$$
It follows that $G=\pi(\lambda D)$ is covered by the $M$ translates $\pi(\lambda\boldsymbol t_i)+A$ of $A$. Translation invariance and subadditivity of $m$ give
$$1=m(G)\leq M\,m(A)\leq24^n m(A).$$
In particular, $m(A)\geq24^{-n}$.

Lemma~\ref{lem:torus-growth} now gives $m(jA)\geq24^{-n/2^{j-1}}$ for every integer $j\geq1$. To make this measure greater than $1/2$, take $j=\lfloor\log_2 n\rfloor+5$. Then $2^{j-1}>8n$, and hence
$$m(jA)\geq24^{-n/2^{j-1}}>24^{-1/8}>\frac12.$$

For every $g\in G$, translation and inversion invariance give $m(g-jA)=m(jA)$, and therefore
$$m\bigl(jA\cap(g-jA)\bigr)\geq 2m(jA)-1>0.$$
A point in this intersection has the form $x=g-y$ with $x,y\in jA$, so $g=x+y\in2jA$. Since $g$ is arbitrary, $2jA=G$.

By convexity of $K$ and additivity of $\pi$, we have $\pi(2j\lambda K)=2jA=G$. Thus $2j\lambda K+\Lambda=\mathbb R^n$, and hence
$$\mu(K,\Lambda)\leq2j\lambda
=\bigl(2\lfloor\log_2 n\rfloor+10\bigr)\mu(D,\Lambda),$$
which proves \eqref{eq:log-upper}.
\end{proof}

\begin{remark}
\label{rem:half-covering}
The implication from $m(jA)>1/2$ to $2j\lambda K+\Lambda=\mathbb R^n$ also follows from Ordentlich, Regev and Weiss \cite[Lemma~2.5]{ORW22}, applied to the convex body $j\lambda K$, since $jA=\pi(j\lambda K)$. Their lemma assumes a strict inequality, but the same conclusion holds when $m(jA)=1/2$ by compactness and connectedness. Put $B=jA$, and fix $g\in G$. The set $B=\pi(j\lambda K)$ is compact. If $B$ and $g-B$ were disjoint, their union would be compact and have measure one. Its complement would be open and have measure zero, and hence empty because every nonempty open subset of $G$ has positive Haar measure. This would express the connected torus $G$ as the union of two disjoint nonempty closed sets, a contradiction. Thus $B\cap(g-B)\neq\varnothing$ for every $g\in G$, giving $2jA=G$ and hence $2j\lambda K+\Lambda=\mathbb R^n$.
\end{remark}

\begin{remark}
\label{rem:covering-constant}
The covering estimates of Rogers and Zong \cite{RZ97} may give sharper versions of \eqref{eq:exponential-cover}. For Proposition~\ref{prop:log-upper}, any bound $N(K-K,K)\leq c^n$ with an absolute constant $c$ suffices for the logarithmic order.
\end{remark}

\section{Other dimensions and a conjecture}
\label{sec:discussion}

Cartesian products allow us to pass from dimensions $2^h$ to arbitrary dimensions. For convex bodies $K_1,K_2$ and full-rank lattices $\Lambda_1,\Lambda_2$ in their respective spaces, we have
\begin{equation}
\label{eq:product-radius}
 \mu(K_1\times K_2,\Lambda_1\times\Lambda_2)
 =\max\{\mu(K_1,\Lambda_1),\mu(K_2,\Lambda_2)\}.
\end{equation}
For every $t>0$,
$$(\Lambda_1\times\Lambda_2)+t(K_1\times K_2)=(\Lambda_1+tK_1)\times(\Lambda_2+tK_2).$$
After translating each body to contain the origin, the covering property of each factor is monotone in $t$. Thus the least scale at which the product covers is the maximum of the two covering radii, proving \eqref{eq:product-radius}.

This identity also shows that $C_{\mathrm{sym}}(n)$ is nondecreasing. For a convex body $K\subset\mathbb R^k$ and a full-rank lattice $\Lambda\subset\mathbb R^k$, where $k<n$, choose $R>0$ so large that
$$\frac1{2R}\leq\mu(K,\Lambda),\qquad \frac1{4R}\leq\mu(K-K,\Lambda).$$
Then $K\times[-R,R]^{n-k}$, with lattice $\Lambda\times\mathbb Z^{n-k}$, has the same covering radius as $K$, and its difference body has the same covering radius as $K-K$. Thus every ratio from a smaller dimension is also realized in dimension $n$.

For $d\geq1$, put $h=\lfloor\log_2d\rfloor$ and $n=2^h$. If $d>n$, take $K_d=P_h\times[-1,1]^{d-n}$, and if $d=n$, take $K_d=P_h$. When $d>n$, we have
$$K_d-K_d=D_h\times[-2,2]^{d-n}.$$
The cube $[-1,1]^{d-n}$ and its difference body have covering radii $1/2$ and $1/4$, respectively, with respect to $\mathbb Z^{d-n}$. By \eqref{eq:exact-radii} and \eqref{eq:product-radius}, we obtain
$$
 \begin{aligned}
 \mu(K_d,\mathbb Z^d)&=\mu(P_h,\mathbb Z^n),\\
 \mu(K_d-K_d,\mathbb Z^d)&=\mu(D_h,\mathbb Z^n).
 \end{aligned}
$$
The same equalities are immediate when $d=n$. Consequently, for $d\geq2$,
$$\frac{C_{\mathrm{sym}}(d)}{\sqrt{\log d}}\geq\frac{\mu(P_h,\mathbb Z^n)}{\sqrt h\,\mu(D_h,\mathbb Z^n)}\sqrt{\frac h{\log d}}.$$
As $d\to\infty$, the first factor on the right tends to $\sqrt{2\pi}$ by \eqref{eq:asymptotic-ratio}. Since $h\leq\log_2d<h+1$, we also have $h/\log d\to1/\log2$. Hence
\begin{equation}
\label{eq:liminf}
 \liminf_{d\to\infty}
 \frac{C_{\mathrm{sym}}(d)}{\sqrt{\log d}}
 \geq\sqrt{\frac{2\pi}{\log2}}.
\end{equation}
The product bodies are generally not simplices. Our explicit integer-simplex examples occur in dimensions $2^h$.

The lower bound in \eqref{eq:liminf} and Proposition~\ref{prop:log-upper} give
$$c\sqrt{\log n}\leq C_{\mathrm{sym}}(n)\leq C\log n\qquad(n\geq2)$$
for absolute constants $c,C>0$. We formulate the following conjecture.

\begin{conjecture}
\label{conj:sqrtlog}
There is an absolute constant $C>0$ such that
$$\mu(K,\Lambda)\leq C\sqrt{\log(2n)}\,\mu(K-K,\Lambda)$$
for every convex body $K\subset\mathbb R^n$ and full-rank lattice $\Lambda\subset\mathbb R^n$.
\end{conjecture}

Together with \eqref{eq:liminf}, Conjecture~\ref{conj:sqrtlog} would give $C_{\mathrm{sym}}(n)=\Theta(\sqrt{\log n})$.
It remains open whether the conjectured upper bound holds for arbitrary simplices.

\section*{Acknowledgements}

During a visit to Berlin, Martin Henk raised the question that led to this work. The authors are grateful to him for subsequent discussions and helpful comments on an earlier version of the paper.

\section*{Declaration of AI use}

The authors used ChatGPT (OpenAI) for language editing, stylistic improvements, and discussions about the presentation of some arguments. ChatGPT was also used to assist in exploring and discussing upper bounds.
All statements and proofs have been verified by the authors, who take full responsibility for the mathematical content.


\begin{thebibliography}{99}

\bibitem{ABG18}
A.~Akopyan, A.~Balitskiy and M.~Grigorev, \emph{On the circle covering theorem by A.~W.~Goodman and R.~E.~Goodman}, Discrete Comput. Geom. \textbf{59} (2018), no.~4, 1001--1009, \doi{10.1007/s00454-017-9883-x}.

\bibitem{ACS26}
D.~Alc\'antara, F.~Criado and F.~Santos, \emph{Covering radii of $3$-zonotopes and the shifted Lonely Runner Conjecture}, Exp. Math. (2026), 1--14, \doi{10.1080/10586458.2026.2651085}.

\bibitem{AH23}
I.~Aliev and M.~Henk, \emph{Minkowski's successive minima in convex and discrete geometry}, Commun. Math. \textbf{31} (2023), no.~2, 35--59, \doi{10.46298/cm.11155}.

\bibitem{BK13}
D.~Berend and A.~Kontorovich, \emph{A sharp estimate of the binomial mean absolute deviation with applications}, Statist. Probab. Lett. \textbf{83} (2013), no.~4, 1254--1259, \doi{10.1016/j.spl.2013.01.023}.

\bibitem{CGdR19}
P.~Candela, D.~Gonz{\'a}lez-S{\'a}nchez and A.~de~Roton, \emph{A {Pl{\"u}nnecke--Ruzsa} inequality in compact abelian groups}, Rev. Mat. Iberoam. \textbf{35} (2019), no.~7, 2169--2186, \doi{10.4171/RMI/1116}.

\bibitem{CFK26}
G.~Codenotti, A.~Freyer and K.~Krivoku\'ca, \emph{Minimal covering bodies and Brunn--Minkowski type inequalities for the covering radius}, preprint, 2026, \href{https://arxiv.org/abs/2606.14442}{arXiv:2606.14442}.

\bibitem{Dad12}
D.~N.~Dadush, \emph{Integer programming, lattice algorithms, and deterministic volume estimation}, Ph.D. thesis, Georgia Institute of Technology, 2012.

\bibitem{Dad19}
\bysame, \emph{On approximating the covering radius and finding dense lattice subspaces}, in Proceedings of the 51st Annual ACM SIGACT Symposium on Theory of Computing (STOC 2019), ACM, New York, 2019, pp.~1021--1026, \doi{10.1145/3313276.3316397}.

\bibitem{GL87}
P.~M.~Gruber and C.~G.~Lekkerkerker, \emph{Geometry of numbers}, second edition, North-Holland Mathematical Library, vol.~37, North-Holland, Amsterdam, 1987.

\bibitem{KL88}
R.~Kannan and L.~Lov{\'a}sz, \emph{Covering minima and lattice-point-free convex bodies}, Ann. of Math. (2) \textbf{128} (1988), no.~3, 577--602, \doi{10.2307/1971436}.

\bibitem{LRR26}
B.~Letwin, V.~Reis and T.~Rothvoss, \emph{Strongly log-concave measures and lattices}, preprint, 2026, \href{https://arxiv.org/abs/2610.29238v2}{arXiv:2610.29238v2}.

\bibitem{ORW22}
O.~Ordentlich, O.~Regev and B.~Weiss, \emph{New bounds on the density of lattice coverings}, J. Amer. Math. Soc. \textbf{35} (2022), no.~1, 295--308, \doi{10.1090/jams/984}.

\bibitem{Rad16}
J.~Radon, \emph{\"Uber eine Erweiterung des Begriffes der konvexen Funktionen mit einer Anwendung auf die Theorie der konvexen K\"orper}, Sitzungsber. Kaiserl. Akad. Wiss. Wien, Math.-Naturwiss. Kl., Abt. IIa \textbf{125} (1916), 241--258.

\bibitem{RR23}
V.~Reis and T.~Rothvoss, \emph{The subspace flatness conjecture and faster integer programming}, in 2023 IEEE 64th Annual Symposium on Foundations of Computer Science (FOCS), IEEE, 2023, pp.~974--988, \doi{10.1109/FOCS57990.2023.00060}.

\bibitem{RS57}
C.~A.~Rogers and G.~C.~Shephard, \emph{The difference body of a convex body}, Arch. Math. (Basel) \textbf{8} (1957), no.~3, 220--233, \doi{10.1007/BF01899997}.

\bibitem{RZ97}
C.~A.~Rogers and C.~Zong, \emph{Covering convex bodies by translates of convex bodies}, Mathematika \textbf{44} (1997), no.~1, 215--218, \doi{10.1112/S0025579300012079}.

\bibitem{Rot26}
T.~Rothvoss, \emph{The subspace flatness conjecture}, in Proceedings of the International Congress of Mathematicians 2026, vol.~6, SIAM, Philadelphia, 2026, pp.~388--406, \doi{10.1137/25M1795546}.

\bibitem{Skr24}
S.~V.~Skresanov, \emph{On directed and undirected diameters of vertex-transitive graphs}, Combinatorica \textbf{44} (2024), no.~6, 1353--1366, \doi{10.1007/s00493-024-00120-4}.

\end{thebibliography}
\end{document}